\documentclass[12pt,reqno]{amsart}
\usepackage{amssymb,latexsym}
\usepackage{graphics,color, amsmath}
\usepackage{graphicx}

\numberwithin{equation}{section}
\newtheorem{theo}{Theorem}
\newtheorem{prop}[theo]{Proposition}

\def\C{\mathbb{C}}
\def\R{\mathbb{R}}

\def\bleu{\textcolor{blue}}
\def\rouge{\textcolor{red}}

\def\Rf{\bleu{R}}
\def\imf{\rouge{I}}
\def\auteur#1{{\sc #1}}
\def\titreref#1{{\em #1}}
\def\vol#1{{\bf #1}}

\title[Combinatorial cell decompositions]{Combinatorial Cellular Decompositions for the Space of Complex Coefficient Polynomials}
\author{Fran\c{c}ois Bergeron}
\address[F. Bergeron]{D\'epartement de Math\'ematiques\\ Universit\'e
  du Qu\'ebec \`a Montr\'eal\\ Montr\'eal, Qu\'ebec, H3C 3P8, CANADA}
\email{bergeron.francois@uqam.ca}
\date{\today}
\thanks{F. Bergeron is supported by NSERC-Canada and FQRNT-Qu\'ebec.}
\begin{document}
\begin{abstract}
 We describe a classification of degree $n$ complex coefficient polynomials with respect to combinatorial patterns that arise from the two real algebraic curves obtained as the zero sets for their real and imaginary part. In particular, we work out explicitly this classification for degree $3$ polynomials, and other special families of polynomials. This work extends to the singular case similar considerations of Martin, Savitt, and Singer for non-singular basketballs.
\end{abstract}
\maketitle

{ \parskip=0pt\footnotesize \tableofcontents}
\parskip=8pt  

\section{Introduction}\label{intro}
As in \cite{martin}, our story starts with Gauss first proof\footnote{Only made completely rigourous much later.} of the fundamental theorem of algebra (see \cite{gauss}). His approach consisted in studying the topology of the two algebraic curves:
\begin{eqnarray}
  \Rf(f)&:=&\{(x,y)\ | \ f(x+y\,i)=i\,t, \ {\rm with}\ t\in \R\}, \qquad {\rm and}\label{bleu}\\
       \imf(f)&:=&\{(x,y)\ |\ f(x+y\,i)=t, \ {\rm with}\ t\in \R\}.\label{rouge}
  \end{eqnarray}
whose intersection is the zero set $Z(f)$ of the polynomial $f$, up to the usual identification between $\R^2$ and $\C$.
In other words, $\Rf(f)$ is the inverse image under $f$ of the imaginary axis, and $\imf(f)$ the inverse image  of the real axis.
We say that the pair $(\Rf(f),\imf(f))$ is an {\em algebraic basketball}, and call $\Rf(f)$  the {\em real component}  of the basketball. Likewise, $\imf(f)$ is the {\em imaginary component} of the basketball. Without loss of generality, we may as well suppose that $f$ is monic of degree $n$.  
Both real algebraic curves $\Rf(f)$ and $\imf(f)$ are harmonic. Thus, by the maximum principle, they have no bounded connected component (``ovals'' in the classical terminology). They are each made up of $n$ branches that go to infinity, and there is no bounded open region of $\C$ whose boundary is entirely made up of portion of the curves $\Rf(f)$ and $\imf(f)$. 

\begin{figure}[ht]
 \begin{center}
\includegraphics[width=35mm]{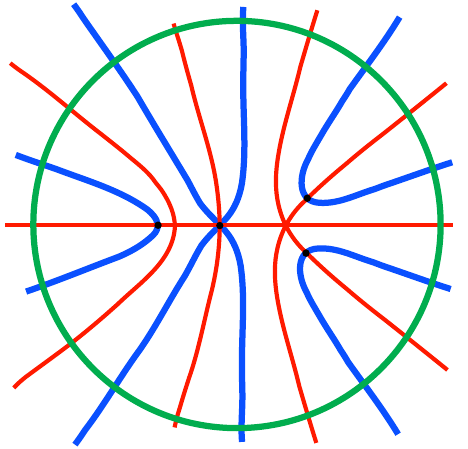}
\end{center}\vskip-10pt
\caption{The curves $\Rf(f)$ and $ \imf(f)$, for 
$f(z)=(z-1)^3(z+1)^2+\frac{2^7\,3^3}{5^5}$.}\label{exemple}
 \end{figure}

\begin{figure}[ht]
 \begin{center}
\includegraphics[height=40mm]{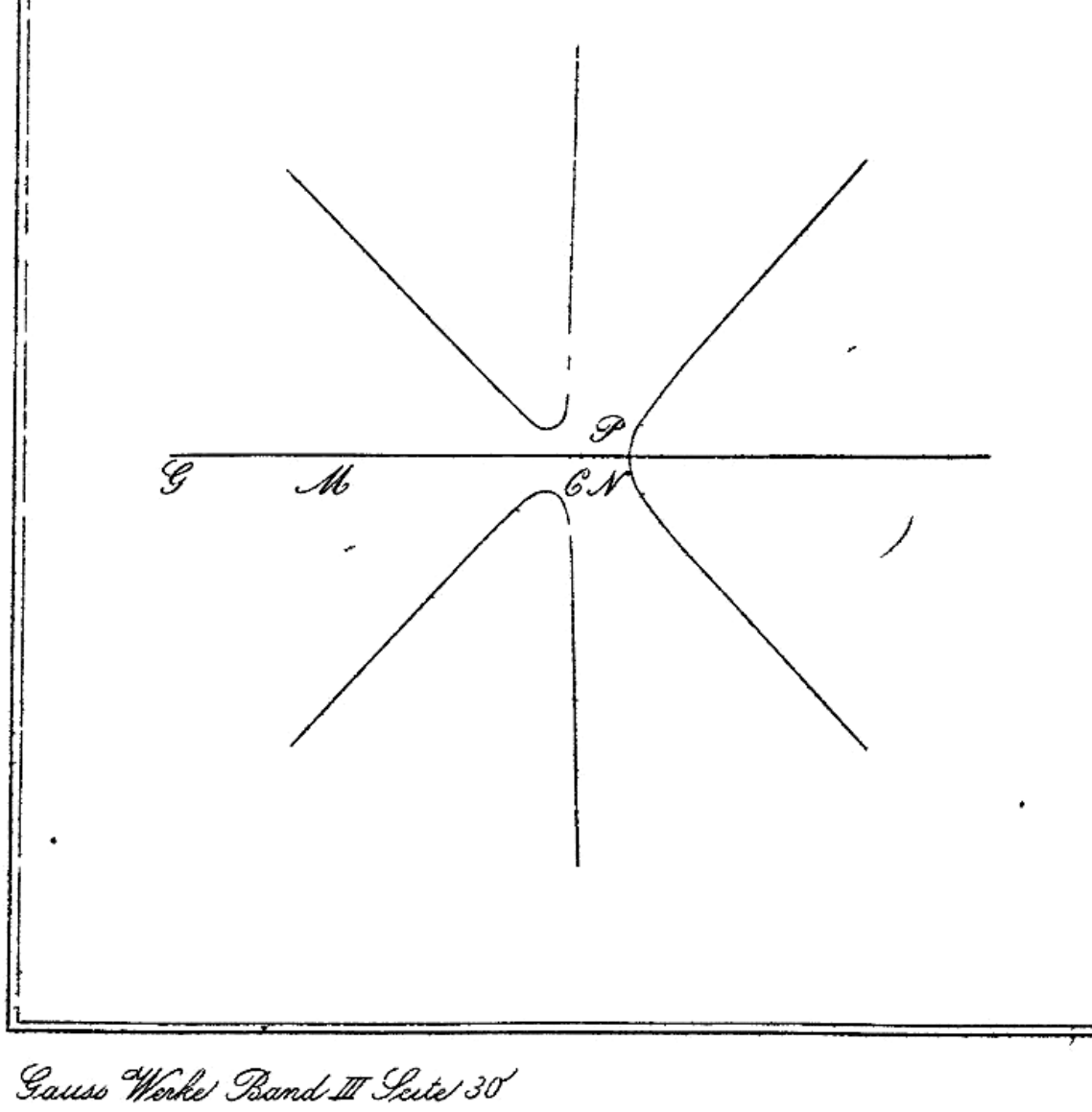}
\end{center}\vskip-15pt
\caption{One of Gauss's figures.}\label{figure_gauss}
 \end{figure}

For any suitably large radius $0$-centered circle $S$, each of the curves $\Rf(f)$ and $\imf(f)$ meet $S$ in $2\,n$ points that are close to the vertices of a regular $2\,n$-gon inscribed in $S$, and all self-intersections or intersections between the two curves lie inside $S$. Moreover, the $2\,n$ points that arise from the intersection $\Rf(f)\cap S$ are interlaced with those coming from $\imf(f)\cap S$. This phenomena is illustrated in Figure~\ref{exemple}, where the thick (blue) curves correspond to $\Rf(f)$ and the thin (red) ones correspond to $\imf(f)$.  A copy of one of of Gauss own figures is presented in Figure~\ref{figure_gauss}. In other words, as the modulus of $z$ goes to infinity, the branches of $\Rf(f)$ (resp. $\imf(f)$) can asymptotically be identified to those of $\Rf(z^n)$ (resp. $\imf(z^n)$), which are half-lines going to infinity. We say that the ``half-branches'' of $\Rf(f)$ become close to those of $\Rf(z^n)$.
Writing $z=re^{i\theta}$ in polar coordinates, the curve $\Rf(z^n)$ is clearly seen to be the set of $n$ lines for which $\cos(n\theta)=0$, so that $\theta=\pi\,(2\,k+1)/(2\,n)$, with $0\leq k\leq n-1$. Similarly, considering the imaginary part, one checks that $\imf(z^n)$ corresponds to the $n$ lines with $\theta=\pi\,k/n$, again with $0\leq k\leq n-1$.

We want to classify polynomials with respect to the ``combinatorial structures'' that can arise as algebraic basketballs $(\Rf(f), \imf(f))$, considering the resulting configuration up to homeomorphism that ``fix'' the $4\,n$ asymptotic directions corresponding to the half-branches of $\Rf(f)$ and  $\imf(f)$. We call the resulting equivalence classes {\em combinatorial basketballs} or simply {\em basketballs}. These combinatorial structures describe the mutual intersections between branches of $\Rf(f)$ (and $\imf(f)$), and how half-branches are comnected.
The branches of $\Rf(f)$ (likewise for $\imf(f)$) are characterized as follows. Starting at some point of intersection between $\Rf(f)$ and $S$, one travels inward along a the chosen branch until one reaches some self-intersection  of $\Rf(f)$. These are going to become internal vertices of our combinatorial structure. Since such self intersection arise at points that are multiple roots of $f(z)-t\,i$, for some $t$ in $\R$, the vertex has an even number of ``edges'' attached to it.  At a vertex one crosses over to continue on the edge that sit opposite to the incoming edge. This is the edge that is separated from the incoming edge by as many edges clockwise as it is counter-clockwise.

Observe that  $\Rf(f)$ (resp. $\imf(f)$) can equivalently be described as the zero set of the real part (resp. imaginary part) of $f(z)$, thus we have
$\Rf(f)=\{z\ | \ {\rm Re}\, f(z)=0\}$, and  $\imf(f)=\{z\ |\ {\rm Im}\, f(z)=0\}$. Moreover, it is clear that
 \begin{eqnarray}
     \Rf(f(z))&=&\Rf(g(z)),\qquad {\rm iff}\qquad f(z)-g(z)\in i\,\R \label{Rcaract}\\
    \imf(f(z))&=&\imf(g(z)),\qquad {\rm iff}\qquad f(z)-g(z)\in \R \label{Icaract}
 \end{eqnarray}
For our classification purpose, may as well suppose that the sum of the roots of $f(z)$ is equal to zero (on top of assuming that $f(z)$ is monic), hence we consider only polynomials of the form
\begin{equation}\label{pol_form}
   f(z)=z^n+c_{2}z^{n-2} +\ldots + c_n.
\end{equation}
 Being interested in the real aspect of the picture, we write the $n-1$ complex coefficients $c_k$ ($k=2,\ldots, n$) in the form $c_k=a_k+b_k\,i$, thus getting $2\,(n-1)$ real  coordinates for our space,  henceforth denoted $\mathcal{P}_n$. In other words, we have identified $\mathcal{P}_n$ with $\R^{2\,(n-1)}$. 

Our objective is to describe a semi-algebraic cell decomposition of $\mathcal{P}_n$, where each cell is characterized by a combinatorial basketball. We also consider a coarser decomposition with cells indexed by ``circular  forests'' (further discussed in the next section) corresponding to the real components of basketballs. This approach makes evident that the set of  basketballs, as well as the set of circular  forests, is endowed with a graded poset structure corresponding naturally to the inclusion of a cell in the closure of another. In particular we consider the enumeration of cells by degree.
Some of the nice results of \cite{savitt} can be considered as the study of certain equivalence classes of closed paths in the space $\mathcal{P}_n$, where these paths only go through cells of co-degree\footnote{Maximal degree minus degree.} $\leq 1$.
  
\section{Combinatorial Basketballs}
\begin{figure}[ht]
 \begin{center}
\includegraphics[width=80mm]{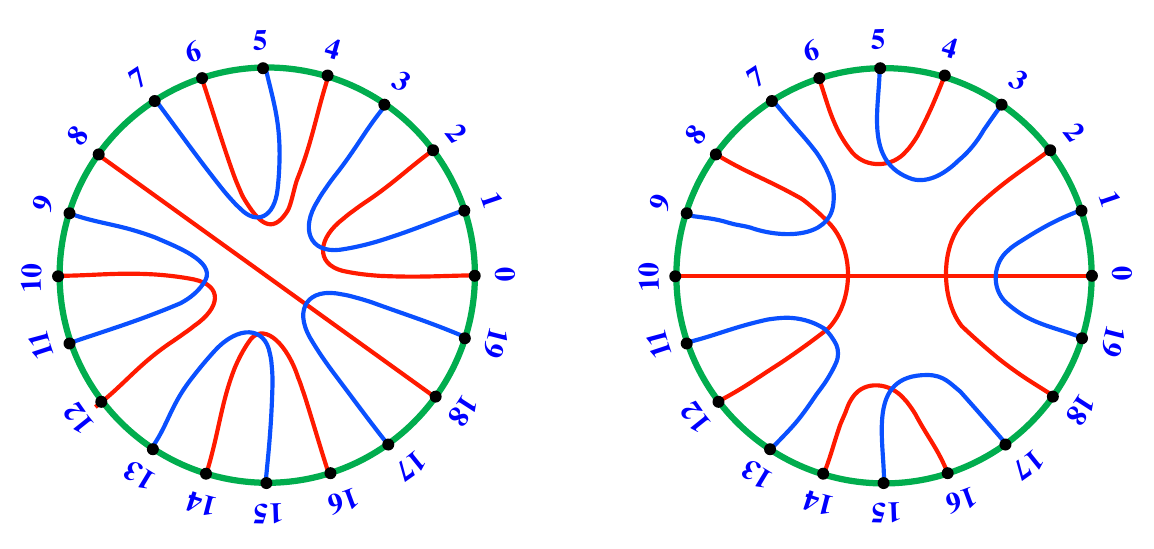}
\end{center}\vskip-10pt
\caption{Nonsingular and singular basketballs.}\label{racines}
 \end{figure}

Our previous discussion implies that both components ($ \Rf(f)$ and $\imf(f)$) of a basketball, are homeomorphic to forests of trees with $2\,n$ leaves (corresponding to the half-branches that go to infinity), with each internal node of even degree. Moreover there is a cyclic structure on the half-edges surrounding any internal vertex.  We thus call them {\em circular  forests}. There is also clearly a unique branch starting at any given leaf, and it must necessarily end at another leaf. In Figure~\ref{racines}, the branch that starts at $0$ ends at $10$. It crosses the branch that goes from $2$ to $18$ and the branch that goes from $8$ to $12$. In this way the leaves of a circular forest are {\em matched}. If all the trees of a circular forest are reduced to two leaves (this is often called a {\em non-crossing} matching) then we say that the forest is {\em non-singular}. Recall that it is a well known fact that the number of non-crossing matching is the Catalan number\footnote{In fact these numbers (formula and generating function) where already known by Euler, see \cite{euler}.}
 \begin{equation}\label{catalan}
    \frac{1}{n+1}\binom{2\,n}{n},
  \end{equation}
 hence it is also the number of non-singular circular forests.

Let us be more precise about the combinatorics of the structures that are involved here. Let $A$ be a finite set of $2\,n$ points on a circle $S$, that we label $0,1,2,\ldots,2\,n$, going counterclockwise starting from the real axis. Consider a partition $\{B_1,B_2,\ldots,B_k\}$ of $A$, into blocks $B_i$.  Two {\em blocks} $B_i$ and $B_j$ are said to {\em cross} if, for some $a<b<c<d$ in $A$, we have
    $$a,c\in B_i,\qquad {\rm and}\qquad b,d\in B_j.$$
In other words, the segments $\overline{ac}$ and $\overline{bd}$ meet inside the disk $D$, for which $S=\partial D$. A partition is said to be {\em non-crossing}, if no two of its blocks cross. A {\em non-crossing matching} of $A$ is a non-crossing partition of $A$ for which all blocks have exactly $2$ elements. For a set $E$ containing $A$, a {\em circular forest} is a graph, $\mathcal{F}=(V,E)$, with all vertex $v\in V$ and edges $e\in E$ drawn inside the disk $D$, in such a way that
  \begin{enumerate}\itemindent=4pt\itemsep=4pt
   \item[(1)]\ each connected component of the graph is a tree whose leaves are the vertex that lie in the set $A:=V\cap S$,
  \item[(2)]\ the partition of $A$, obtained by choosing blocks to be the set of leaves of connected components, is non-crossing,
  \item[(3)]\ each inside vertex, namely those  in $V\setminus A$, has even degree.
 \end{enumerate}
We say that $S$ is the {\em bounding} circle of the circular forest, and condition (1) says that the leaves all lie on the bounding circle.
In view of our previous discussion, such a circular forest is made up of $n$ branches. These branches establish a matching (we sometines say a $n$-matching) between the $2\, n $ leaves of the forest. It is not hard to check that this matching entirely characterizes the circular forest, and branches can be identified with pairs $\{a,b\}$ of elements of $A$. We sometimes denote such pair simply as $ab$.

With our previous terminology, two branches $ab$ and $cd$ cross if 
    $$a<c<b<d.$$ 
  In other terms, they share an internal vertex. More generally, for $1\leq k\leq m$,  several branches $a_kb_k$, of a forest $F$,  go through one common vertex if and only if we have 
\begin{equation}\label{sommet}
   a_1<a_2<\ldots <a_m<b_1<b_2<\ldots <b_m,
  \end{equation}
 for some correct choice of the labels $a_k$ and $b_k$. Thus the vertices of a forest $F$ are characterized as maximal sets of branches 
$\mathcal{B}=\{a_1b_1,a_2b_2,\ldots, a_mb_m\}$ such  that condition (\ref{sommet}) hold. 
To characterize the $n$-matchings of the $2\,n$ points of $A$  associated to circular forests, we proceed as follows. Let us say that a set of pairs of $A$
     $$\mathcal{B}=\{a_1b_1,a_2b_2,\ldots, a_mb_m\},$$
 is a {\em crossing cycle} if $a_kb_k$ crosses $a_{k+1}b_{k+1}$ for all $k$, with the convention that $a_{m+1}=a_1$ and $b_{m+1}=b_1$. 
 The matchings  associated to circular forests are exactly those  that contain only crossing cycles that
 correspond to vertices, i.e.: crossing cycles satisfying condition (\ref{sommet}).
Another notion along these lies that will be useful later is that of  {\em crossing set} $\mathcal{C}(F)$ of a circular forest $F$. This is the set
 of pairs $\{j,k\}$, of leaves of $F$, such that the branch on which $j$ lies crosses the branch on which $k$ lies. 
 
Two circular forests with $2\, n$ leaves, $\Rf$ and $\imf$, are said to be {\em orthogonal} if and only if 
\begin{enumerate}\itemindent=4pt\itemsep=4pt
 \item[(a)]  their respective leaves are interlaced around the same bounding circle,  
 \item[(b)] each branch of $\Rf$ meets at least one branch of $\imf$, and vice-versa.
  \item[(c)] the union of $\Rf$ and $\imf$ is a circular forest. 
 \end{enumerate}
 A {\em combinatorial $n$-basketball} is an {\em orthogonal} pair $(\Rf,\imf)$ of  circular forests, each with $2\,n$ leaves. We intend combinatorial basketballs to be representatives of homotopy classes of algebraic basketballs in the sense of section~\ref{intro}. A combinatorial basketball  $(\Rf,\imf)$ is said to be {\em non-singular} if and only if both $\Rf$ and $\imf$ are non-singular.
 As discussed in \cite{martin}, the number of non-singular combinatorial basketballs\footnote{In their paper called simply basketballs, since the singular cases are not considered.}  is
 \begin{equation}\label{formule_non_sing}
    \frac{1}{3\,n+1}\binom{4\,n}{n}
  \end{equation}
Observe that for two polynomials $f(z)$ and $g(z)$ that differ by some real factor homothety, i.e.:
$g(z)=t^n f(z/t)$, with $t>0$, one gets the same class of algebraic basketball.
Thus, up to a rescaling, we may as well suppose that $S$ is the unit circle and that $A$ is the set of points
     $$\{\exp((2k+1)\,\pi\,i/2n)\ |\ 0\leq n-1\ \}$$
in the case of $\Rf(f)$, and 
        $$\{\exp(k\, \pi\,i/n)\ |\ 0\leq n-1\ \}$$
in the case of $\imf(f)$. As before, these points are labelled form $0$ to $4n-1$ going counterclockwise. Each half-branch (corresponding to one of these points) is connected through a branch to a unique other half-branch. We get in this way a natural decomposition of $\Rf$ (and $\imf$) into $n$ branches. The fact that in an algebraic basketball the union of $\Rf(f)$ and $\imf(f)$ is a circular forest is readily deduced from the simple observation that
\begin{equation}\label{union}
   \imf(f^2)=\Rf(f)\cup \imf(f).
 \end{equation}
 
\section{Cellular decompositions}
We are interested in describing explicit cellular decompositions $\Gamma=\Gamma_{\mathcal{Q}}$ of particular subspaces $\mathcal{Q}$ of the space $\mathcal{P}_n$:
     $$\mathcal{Q}=\sum_{\gamma\in \Gamma} \gamma,$$
  where we use summation to underline that we have a disjoint union of cells. 
 Cells of these (finite) decompositions correspond to classes of algebraic basketballs, and these are naturally indexed by (combinatorial) basketballs.  Each cell $\gamma$ is semi-algebraically homeomorphic to some open hypercube of dimension $d=d(\gamma)$, and we say that $\gamma$ is a $d$-cell. We can thus enumerate the cells by this {\em degree} $d$. The resulting {\em cell enumerator polynomial}  is defined to be
\begin{equation}\label{cell_enumeration}
   \Gamma(x):=\sum_{\gamma\in \Gamma} x^{d(\gamma)}.
\end{equation}
We may very well consider different cell decompositions of the same space. Indeed, we will be interested in both the cell decomposition corresponding to classifying degree $n$ polynomials only with respect to the type of the real component of their basketball, and the cell decomposition that takes into account the whole basketball. The first decomposition is denoted $\mathcal{F}_n$, and called the {\em forest-cell-decomposition}, since its cells are identified with circular forests. The second decomposition,  denoted $\mathcal{B}_n$, is  the {\em basketballs-cell-decomposition}.
Thus, the coefficient of $x^{2\,(n-1)}$ in  $\mathcal{F}_n(x)$ is the Catalan number, whereas in  $\mathcal{B}_n(x)$ it is number given by formula (\ref{formule_non_sing}). 

There is a natural graded poset structure on the cells of a decomposition. Namely one sets $\gamma_1\leq \gamma_2$ if $\gamma_1$ is included  in the closure of $\gamma_2$, and the grading is given by the degree of the cells. In terms of circular forests, the order can be characterized as follows. The cell corresponding to a circular forest $F_1$ is smaller than the cell associated to $F_2$, if and only if the crossing set of $F_1$ contains the crossing set of $F_2$. In formula
\begin{equation}
  F_1\leq F_2\qquad \mathrm{iff}\qquad \mathcal{C}(F_2)\subseteq \mathcal{C}(F_1).
\end{equation}
The minimal cell (of degree $0$) is the class of the polynomial: $z^n$. For the basketball decomposition, the order is just extended component-wise. For each cell $\gamma$, in the decomposition considered, the closure of $\gamma$ is the union of the cells that are smaller than $\gamma$, i.e.:
 \begin{equation}\label{fermeture}
     \overline{\gamma}=\sum_{\gamma'\leq \gamma} \gamma'.
  \end{equation}
 
\section{Non singular basketballs}
\begin{figure}[ht]
 \begin{center}
\includegraphics[width=90mm]{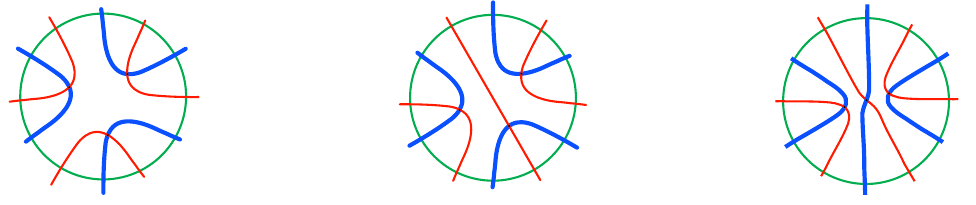}
\end{center}\vskip-10pt
\begin{picture}(0,0)(0,0)
\small
\put(-150,0){(A.1)\ $z^3+(1-i)$} \put(-50,0){(A.2)\ $z^3-3\,e^{i\,\pi/3}z+2$}  \put(70,0){(B)\ $z^3-3\,(1+i)\, z$}    
\end{picture}
\caption{Types of non-singular basketballs for $n=3$.}\label{vingtdeux}
 \end{figure}
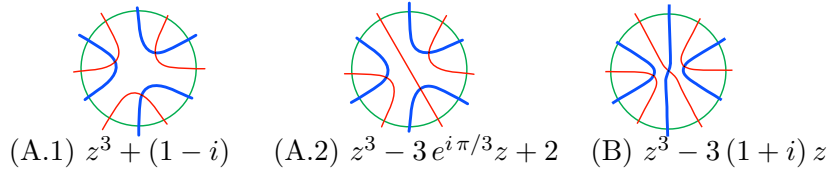

  In generic cases, both circular forests $\Rf(f)$ and $\imf(f)$ are non-singular, and they are maximal elements of the poset considered above. More precisely, in a dense open subset of  $\mathcal{P}_n$, the circular forest  $\Rf(f)$ (resp. $\imf(f)$) is a  non-crossing matching  between the $2\,n$ points 
    $$\{1,3,5,\ldots, 4n-1\},\qquad {\rm (resp. }\   \{0,2,4,\ldots, 4n-2\}).$$
 In other words, both $\Rf(f)$ and  $\imf(f)$ consist of $n$ non-intersecting branches and, by our previous discussion, each branch of $\Rf(f)$ crosses one and only one branch of $\imf(f)$. It is a theorem of \cite{martin} (called the ``Inverse Basketball Theorem'') that each possible non-singular combinatorial basketball can actually be realized (up to homotopy) as a (non-singular) algebraic basketball. Hence non-singular combinatorial basketballs are natural representatives of our classes of algebraic basketballs. In Figure~\ref{vingtdeux} we give examples of all ``types'' of non-singular basketballs for $n=3$. 
Observe that, beside properties (\ref{Rcaract}) and (\ref{Icaract}), there are interesting symmetries of the space $\mathcal{P}_n$. Indeed, for $f(z)$ in $\mathcal{P}_n$, we have
     $$i\,f(e^{-i\,\pi /2n}\,z) \in \mathcal{P}_n,$$
 although this symmetry exchanges the role of $\Rf(f)$ and $\imf(f)$.
More precisely, let $\psi_n$ be the rotation of a basketball by an angle of $\pi/2\,n$, thus relabelling the leaves 
    $$k\mapsto (k+1\mod 4\,n).$$  
  We observe that
\begin{eqnarray}
    \Rf(i\,f(e^{-i\,\pi /2\,n}\,z)&=&\psi_n(\imf(f)),\\
    \imf(i\,f(e^{-i\,\pi /2\,n}\,z)&=&\psi_n( \Rf(f)).
  \end{eqnarray}
Considering the further symmetry
  \begin{equation}\label{reflexion}
      \Rf(\overline{ f(\overline{z})})=\tau(\Rf(f))
 \end{equation}  
 with $\tau$ being the reflection that sends $k$ to $4\,n-1-k$, we get an explicit action of the dihedral group (of order $4\,n$) on polynomials which respects the basketball cell decomposition. We can thus apply the classical techniques of orbit enumeration to count basketballs. 
To illustrate, from the polynomials of Figure~\ref{vingtdeux} we get 22 non singular basketballs adding up the relevant orbits, which are of respective size $4$, $6$ and $12$. These are readily seen to be the only possibilities, and the total count does correspond to formula~(\ref{formule_non_sing}).

\section{Singular basketball locus}\label{sec_sing}
In order to give a complete description of all possible algebraic basketballs, we need to extend our classification to those that exhibit singularity, i.e.:  either $\Rf(f)$ or $\imf(f)$ is singular. We first observe  that $\Rf(f)$ is singular if and only if the following {\em discriminant}\footnote{This discriminant with respect to the variable $z$ is a polynomial in $t$, normalized to be monic. Equations~(\ref{Lyaschko}) and (\ref{Looinjenga}) involve the  Lyaschko-Looinjenga map 
    $$LL:f(z)\mapsto {\rm Disc}_z(f(z)-t).$$
As further discussed in \cite{zvonkin}, this is a degree $n^{n-2}$ mapping.} vanishes
\begin{equation}\label{Lyaschko}
   {\rm Disc}_z(f(z)-i\,t)=0
\end{equation}
for some real value of $t$. Similarly $\imf(f)$ is singular if and only if
\begin{equation}\label{Looinjenga}
  {\rm Disc}_z(f(z)-t)=0,
\end{equation}
for some $t\in\R$. To characterize in the space $\R^{2\,(n-1)}$ the locus of polynomials 
     $$f(z)=z^n+(a_2+b_2 \,i)\,z^{n-2} + \ldots + (a_{n-1}+b_{n-1}\,i)\, z + (a_{n}+b_{n}\,i)$$
having a singular $\Rf(f)$, one eliminates the variable $t$ in the two equations that arise from considering the real part and the imaginary part of equation~(\ref{Lyaschko}). This process gives a polynomial in the variables $a_k$ and $b_k$, which is denoted $\Delta_n=\Delta_n(f)$, and $\Rf(f)$ is singular if and only if $\Delta_n(f)=0$. The analogous polynomial for $\imf(f)$ is denoted $\Delta'_n$.
 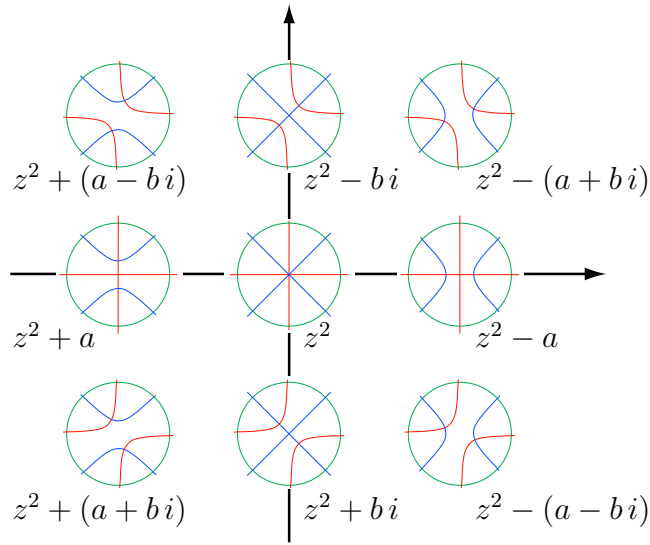
\begin{figure}[ht]
 \begin{center}
\includegraphics[width=80mm]{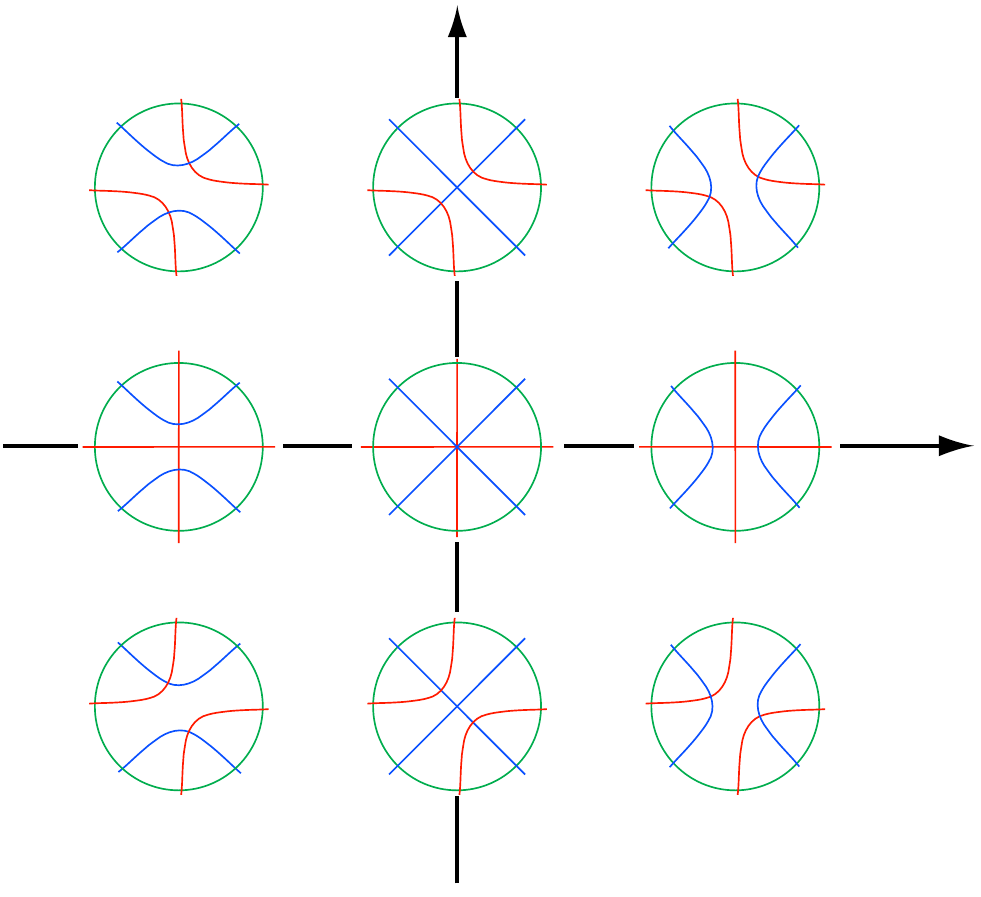}
\begin{picture}(0,0)(120,18)
\put(0,30){$z^2+b\,i$}
\put(0,95){$z^2$}
\put(0,155){$z^2-b\,i$}
\put(65,30){$z^2-(a-b\,i)$}
\put(65,95){$z^2-a$}
\put(65,155){$z^2-(a+b\,i)$}
\put(-110,30){$z^2+(a+b\,i)$}
\put(-110,95){$z^2+a$}
\put(-110,155){$z^2+(a-b\,i)$}
\end{picture}\end{center}
\vskip-10pt
\caption{Basketballs for $n=2$.}\label{sing2}
 \end{figure}

Let us first illustrate this with the space $\mathcal{P}_2$ of polynomials $f(z)=z^2-(a+b\,i)$ of degree $2$, identified with points $(a,b)$ in $\R^2$. 
The only such polynomials that give rise to singular basketballs are those of the form
   $$z^2-b\,i,\qquad {\rm or}\qquad z^2-a,\qquad (b\in\R).$$
Indeed,  one calculates that $\mathrm{Disc}(z^2-(a+b\,i)-i\,t)=t+(b-a\,i)$, and thus $\Delta_2=a$. Hence the singular locus for $\Rf(f)$ is  the $b$-axis in the $(a,b)$-plane. Likewise, one calculates that $\Delta'_2=b$ so that the singularity locus for $\imf(f)$ is the $a$-axis.  The corresponding cellular decomposition of $\mathcal{P}_2$ is thus  as follows:
\begin{enumerate}\itemindent=4pt\itemsep=4pt
\item[(a)] one $0$-cell made up of the polynomial $z^2$, which lies in both singularity locus;
\item[(b)] four $1$-cells respectively containing the polynomials $z^2+b\,i$, $z^2-b\,i$, $z^2+a$ and $z^2-a$, assuming that $a>0$ and $b>0$;
\item[(c)] four $2$-cells respectively containing the polynomials $z^2+(a+b\,i)$, $z^2+(a-b\,i)$, $z^2-(a+b\,i)$ and $z^2-(a-b\,i)$, still assuming that $a>0$ and $b>0$. These cells correspond to  the non-singular situations.
\end{enumerate}
We conclude that the cell enumerator polynomial $\mathcal{B}_2$ is
   $$\mathcal{B}_2=1+4\,x+4\,x^2.$$

\section{Critical values}
We can generalize the analysis of the example of the previous section to study the basketball cell decomposition of the space of polynomials of the form $f(z)-(a+b\,i)$ for any fixed polynomial $f(z)$. Indeed, the behaviour of the basketball of the polynomial $ f(z)-\alpha$, as $\alpha=a+b\,i$ varies in $\C$, can be entirely described in terms of the relative positions of the ``critical values'' of $f$. Recall that  the {\em critical values} of $f$ are the complex numbers $\alpha=f(\omega)$, such that $\omega$ is a root of the derivative of $f$. In other terms,  critical values of $f$ are the $\alpha=f(\omega)$ such that the polynomial $f(z)-\alpha$ has at least one multiple root, since both $f(z)-\alpha$ and its derivative vanish at $\omega$. Observe that the self intersections of the real component $\Rf(f)$, as well as those of the imaginary component $\imf(f)$, occur at roots of $f'(z)$. It follows that the total number (with multiplicities correctly counted) of these self-intersection is at most $n-1$.

\begin{figure}[ht]
 \begin{center}
\includegraphics[width=80mm]{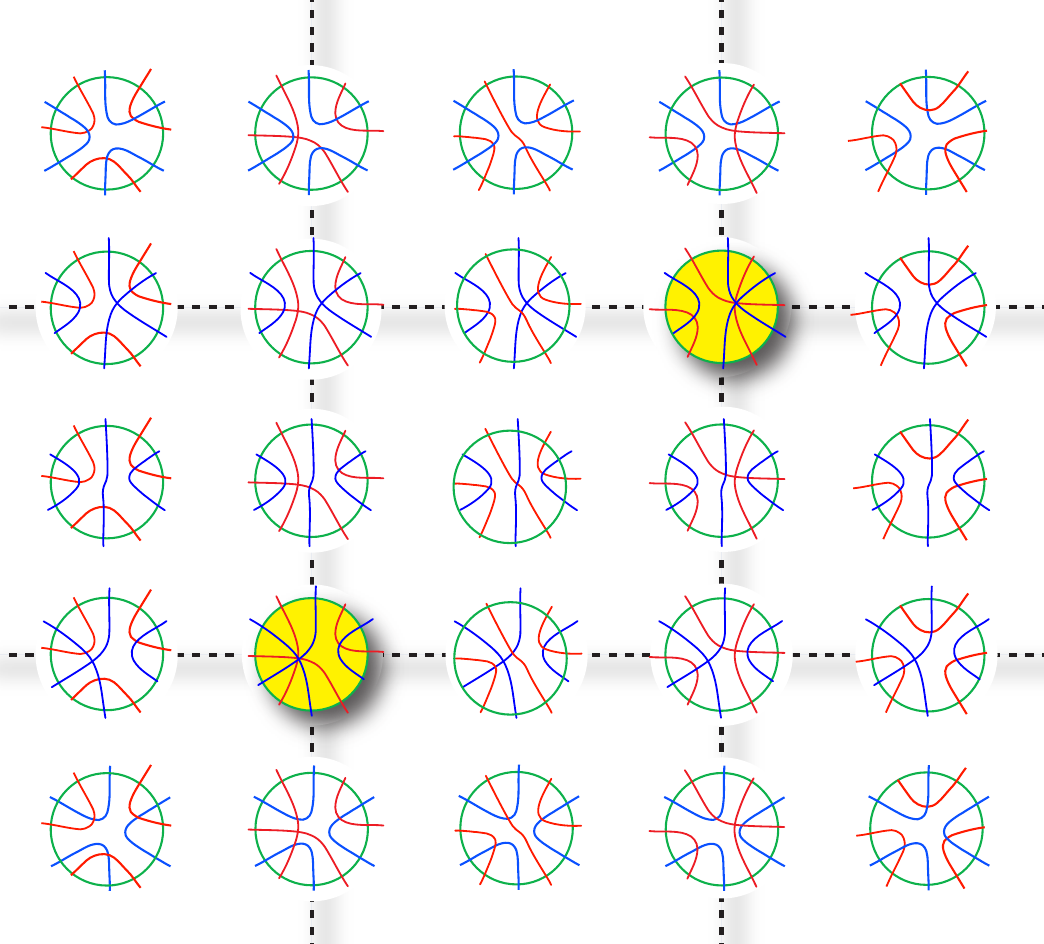}
\end{center}
\vskip-10pt
\caption{Cellular decomposition for $f(z)-(a+b\,i)$.}\label{generique3}
 \end{figure}

\begin{figure}[ht] 
 \begin{center}
\includegraphics[width=90mm]{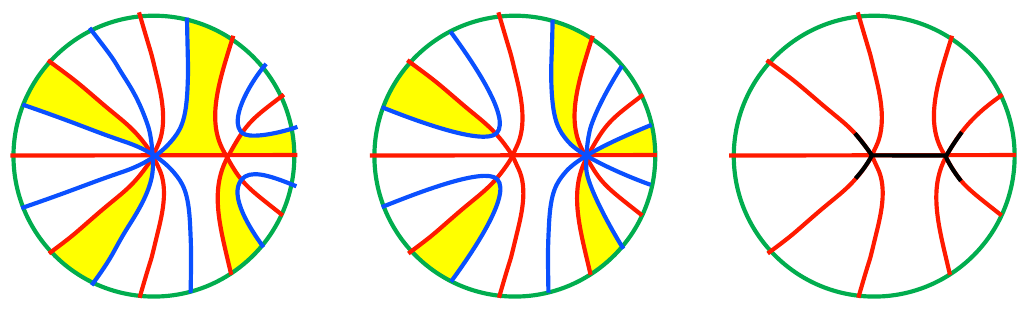}
\end{center}
\vskip-10pt
\caption{Basketballs of $f(z)$ and $f(z)-1$, and resulting dessin d'enfant.}\label{fig_shabbat}
 \end{figure}

At a critical value $\alpha=a+b\,i$, the basketball of $f(z)-\alpha$ is singular, both in its real and imaginary component. Moreover, along the vertical line (resp. horizontal line) $\beta=a+t\,i$ (resp. $\beta=t+b\,i$), for $t$ in $\R$, the real component (resp. imaginary component) of the basketball of $f(z)-\beta$ remains singular in its real (resp. imaginary) component. In fact, this is the whole story. More precisely, let $a_1<a_2< \ldots< a_\ell$  be the distinct values that appear as real parts of critical values of $f(z)$. 
Likewise, let $b_1<b_2<\ldots < b_m$ be the distinct values that appear as imaginary parts of critical values of $f(z)$
With $0\leq j\leq m$ and $0\leq k\leq \ell$, the cells of the basketball decomposition of the $(a,b)$-plane (i.e.: the space of polynomials of the form $f(z)-(a+b\,i)$) are as follows\footnote{To simplify our presentation, we assume that $a_0=b_0=-\infty$ and $a_{\ell+1}=b_{m+1}=+\infty$.}:
\begin{enumerate}\itemindent=4pt\itemsep=4pt
   \item[(1)] restricting both $j$ and $k$ to be larger or equal to $1$, each of the points $(a_j,b_k)$ is a cell of degree 0;
   \item[(2)] each of the open segments going from $(a_j,b_k)$ to $(a_j,b_{k+1})$, and from $(a_j,b_k)$ to $(a_{j+1},b_k)$, are cells of degree $1$; and
   \item[(3)] the rectangles having south-west corner $(a_j,b_k)$ and north-east corner $(a_{j+1},b_{k+1})$ are the degree $2$ cells.
\end{enumerate}
 The cell enumerator for the $(a,b)$-plane under study is thus
\begin{equation}
   m\,\ell + (2\,m\,\ell+m+\ell)\,x+ (m\,\ell+m+\ell+1)\,x^2.
  \end{equation}
For instance, one calculates that the critical values of the polynomial $f(z)=z^3+6\,i\,z$ are $-4-4\,i$ and $4+4\,i$.
The corresponding basketball cell decomposition of the $(a,b)$-plane  is illustrated in Figure~\ref{generique3}, with the two critical values highlighted. 
From its characterization as the inverse image of the imaginary axis, it is easy to check that the following holds for $\Rf(f)$.
 \begin{prop}
For a polynomial $f$, the circular forest $\Rf(f)$ {\rm (}resp. $\imf(f)${\rm )}  is a tree if and only if the critical values of $f$ all have  real part {\rm (}resp. imaginary part{\rm )} equal to 0.
 \end{prop}
 This links our discussion to the notion of {\em dessins d'enfants}\footnote{It seems that this notion was already present in Klein's work (see \cite{klein}) under the name {\em linienzuges}, as observed by le Bruyn (see \cite{lebruyn}).} of Grothendieck (see \cite{grothendieck}), more particularly to the study of ``Shabat polynomials''. Recall that $f$ is a {\em Shabat polynomial} if and only if it has precisely two critical values. Usually it is assumed that these critical values are $0$ and $1$, and one is interested in the combinatorial structure obtained as the inverse image  of the interval $[0,1]$ under $f$. This is the associated {\em dessin d'enfant}, and it is a tree-shaped subset of the imaginary component $\imf(f)$. Indeed, the dessin d'enfant of $f(z)$ is the portion of the imaginary component that corresponds to paths followed by zeros of $f(z)-t$ as $t$ goes from $0$ to $1$.  
 To illustrate this last assertion one can calculate that, for the polynomial
 $$f(z)=z^3\,(z^2-(5/2)\,6^{1/5}\,z+(5/3)\,6^{2/5}),$$
the basketballs of the Shabat polynomial $f(z)$ and $f(z)-1$ are as given in the first two circles of Figure~\ref{fig_shabbat}, 
and the dessin d'enfant of $f(z)$ is outlined in black in the third. 

\begin{figure}[ht] 
 \begin{center}
\includegraphics[width=30mm]{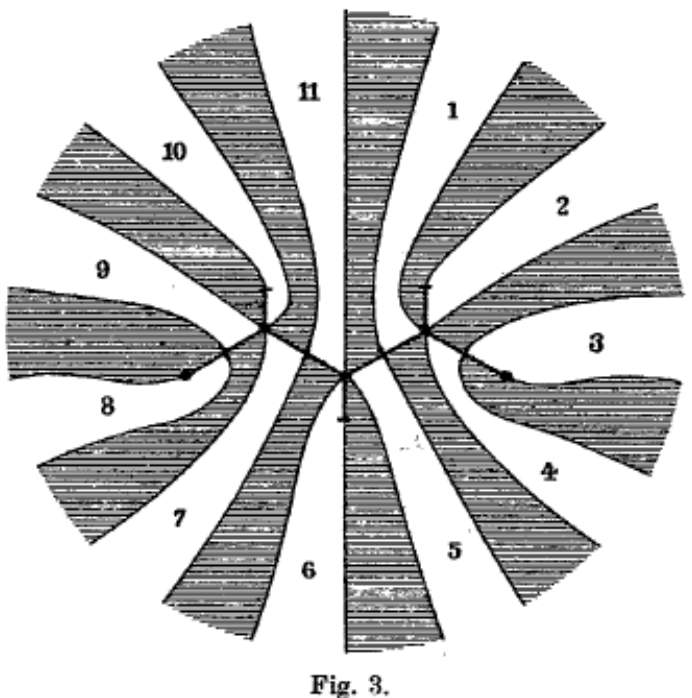}
\end{center}
\vskip-10pt
\caption{One of Klein's illustrations.}\label{fig_klein}
 \end{figure}

To make clearer the combinatorial link between basketballs and dessin d'enfants,  we orient the branches $\{a,b\}$ of basketballs from $a$ to $b$, if $b\equiv 0,1\ ({\rm mod}\ 4)$. This makes sense since there is one and only one of the two extremities of a branch satisfying this property (see \cite{martin}). In the imaginary (resp. real) component case, traveling on a branch in this positive direction corresponds to inverse images of values going from $-\infty$ to $\infty$ along the real (resp. imaginary) axis. 
Observe that, for any $k\equiv 0 \ ({\rm mod}\ 4)$, the internal region that corresponds to the arc on  the circle $S$ going  from $k$ to $k+1$ is
the inverse image of the quadrant $\R^+\times \R^+$. These are the shaded regions in Figure~\ref{fig_shabbat}. Klein used exactly this kind of illustrations (see Figure~\ref{fig_klein} for an example taken from \cite{klein})  to make clearer the action of  monodromy groups (in modern terminology).

\section{Degree 3 polynomials: the real component}\label{sec_real}
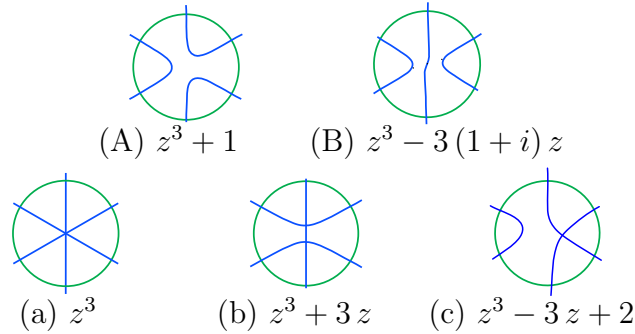
\begin{figure}[ht]
 \begin{center}
\includegraphics[width=80mm]{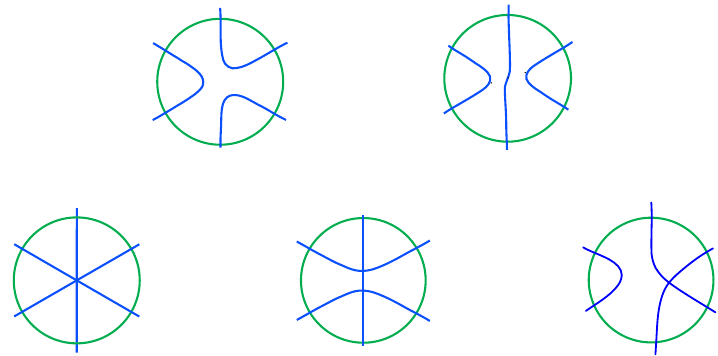}
\begin{picture}(0,0)(225,0)
 \put(30,55){(A)\  $z^3+1$}
  \put(110,55){(B)\ $z^3-3\,(1+i)\,z$}
 \put(0,-10){(a)\ $z^3$}
  \put(75,-10){(b)\ $z^3+3\,z$}
  \put(155,-10){(c)\ $z^3-3\,z+2$}
\end{picture}
\end{center}
\caption{Circular forest representatives for $n=3$.}\label{fig_forest_trois}
 \end{figure}

\begin{figure}[ht]
 \begin{center}
\includegraphics[width=40mm]{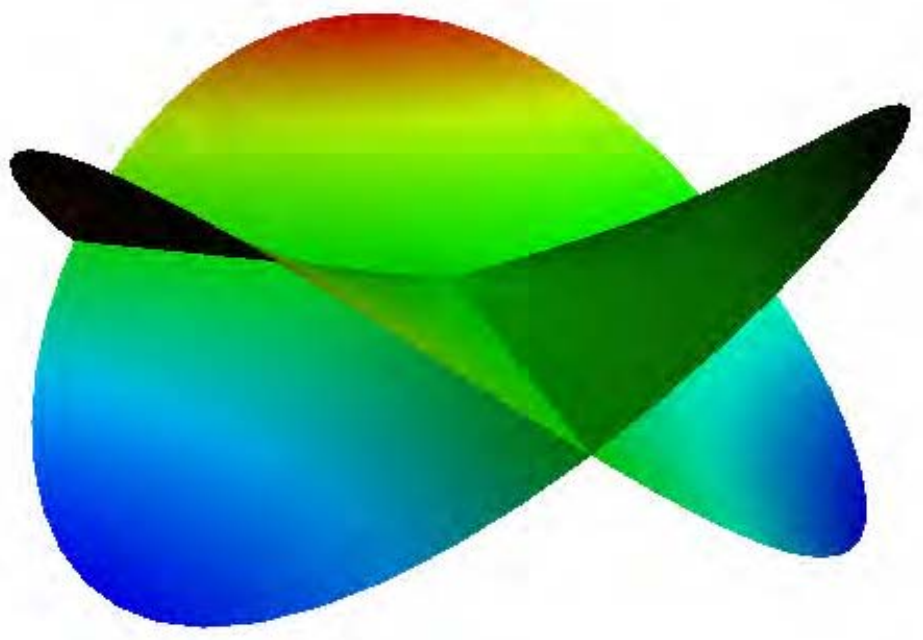}
\end{center}\vskip-10pt
\caption{Real component singularity locus $Z(\Delta_3)$.}\label{fig_sing_trois}
 \end{figure}

Before going on with more general results, let us study  the situation for $n=3$. 
We first consider only the classification of polynomials 
     $$f(z)=z^3-3\,(a+b\,i)z +2\,(c+d\,i)$$ 
  with respect to the type of its real component $\Rf(f)$. The integers multiplying the coefficients make upcoming calculations simpler.  Under the action of the pertinent dihedral group, the corresponding set of cells  breaks up into five orbits for which we give representatives\footnote{The picture in (b) is the reason for the name ``basketballs''.} in Figure~\ref{fig_forest_trois}.   
By the process outlined above, one calculates that $\Rf(f)$ is singular if and only if one chooses $a$, $b$ and $c$ such that $\Delta_3=0$, with
\begin{equation}\label{condition_trois}
      \Delta_3:=(2\,c^2-a\,(a^2-3\,b^2))^2-(a^2+b^2)^3
\end{equation}
Observe that $\Delta_3$ only depends on $a$, $b$ and $c$. Thus the forest-cell-decomposition of  $\mathcal{P}_3$ is independent of $d$, and can essentially be described in $\R^3$.  
Observe also that, setting $a=r\,\cos(\theta)$ and $b=r\,\sin(\theta)$, $\Delta_3$ factors as
      $$ \Delta_3=4\,(c^2-r^3\cos^2(3\,\theta/2))\,(c^2+r^3\sin^2(3\,\theta/2)).$$
 Since the second factor is non-negative for all real $c$ and $r>0$, the condition that $\Delta_3\geq 0$ is equivalent to 
  \begin{equation}\label{cond_cyl}
    c^2\geq r^3 \cos^2(3\, \theta/2).
    \end{equation}
In particular the $\Rf$ singularity locus can be parametrized as follows
     $$Z(\Delta_3)=\left\{\ (r\cos(\theta) ,r\sin(\theta) , r^{3/2}\cos(3\,\theta/2),d)\ |\ r\geq 0,\ 0\leq \theta< 4\,\pi\right\}.$$
The restriction of this surface to $\R^3$ (setting $d=0$) is illustrated in Figure~\ref{fig_sing_trois}, with the $c$ axis pointing up.
 It separates the $(a,b,c,d)$-space $\R^4$ into $5$ open semi-algebraic regions ($4$-cells), each of which characterized by a non-crossing matching $\Rf$ which is non-singular. The sign of  $\Delta_3$ is clearly constant in these regions and it changes when one goes through a wall\footnote{A co-degree $1$ cell.} of the surface $\Delta_3=0$. Going through a point of self-intersection of $\Delta_3=0$ may correspond to crossing several walls. To effectively calculate the matching associated to a cell, we clearly need only do it for one of the polynomials in this cell. Moreover, since the action of the relevant dihedral group clearly translates into an action of this group on matchings, we actually need only calculate the matching for one of the cells in a given orbit. The $4$-cells can now be characterized as follows:
 \begin{enumerate}\itemindent=4pt\itemsep=4pt
 \item[(A)] We begin with the two $4$-cells that correspond to cases for which we have either $c> r^{3/2}\cos(3\,\theta/2)>0$, or 
 $c< r^{3/2}\cos(3\,\theta/2)<0$.   For instance, the first cell contains the polynomial $z^3+1$ (see Figure~\ref{fig_forest_trois} (A)).
 The two cells are bijectively mapped onto one another by the correspondance $c\mapsto -c$, and they correspond respectively to the non-crossing matchings:
   $$\{\bleu{\{1,3\},\{5,7\},\{9,11\}}\},\qquad \mathrm{and}\qquad \{\bleu{\{1,11\},\{3,5\},\{7,9\}}\},$$
 the first of which being the non-crossing matching characterizing the cell containing $z^3+1$. It is easy to check that the matching above indeed corresponds to how the
branches of the real component curve $x(x^2-3\,y^2)=0$ connect the $6$ points on the circle.
 \item[(B)] Then there are three other $4$-cells, all similar, for which we have 
\begin{equation}\label{cond_bouche}
    -r^{3/2} \cos(3\, \theta/2)<\ c\ <r^{3/2} \cos(3\, \theta/2),
 \end{equation}
and with $\theta$ varying in one of the following intervals: \\[-10pt]
  \begin{enumerate}\itemindent=15pt\itemsep=4pt
     \item[(1)] $-\pi/3<\theta<\pi/3$, for the first region (see Figure~\ref{fig_forest_trois}-(B));
     \item[(2)] $\pi<\theta <5\, \pi/3$, for the second region; and
     \item[(3)] $7\,\pi/3<\theta<3\,\pi$, for the third region. \\[-10pt]
    \end{enumerate} 
  Observe that in each of these regions $\cos(3\, \theta/2)$ is positive, so that condition (\ref{cond_bouche}) makes sense.
  These cells are clearly mapped to one another by a rotation of $2\pi/3$ around the $c$-axis. They respectively are characterized by the non-crossing matchings
   $$\{\bleu{\{1,11\},\{3,9\},\{5,7\}}\},\quad\{\bleu{\{1,7\},\{3,5\},\{9,11\}}\},\quad \mathrm{and}\quad \{\bleu{\{1,3\},\{5,11\},\{7,9\}}\}.$$  
  as seen by calculating the branches of the real component curve  
     $$x\,(x^2-3\,y^2-1)=0$$ 
  corresponding to the polynomial $z^3-z$ which lies in the first cell.   
 \end{enumerate} 
The whole cellular decomposition of $\mathcal{P}_3$ consists of a total of $15$ cells, with the remaining $10$ cells arising from a natural decomposition of the surface $Z(\Delta_3)$ itself, which goes as follows. We start with the decomposition of  the locus of self-intersections of the surface. Considering $\Delta_3$ as a polynomial in $c$, this locus corresponds to values of $a$ and $b$ for which there $\Delta_3$ has multiple roots. Hence the discriminant of this polynomial in $c$ must vanish, so that we find that
\begin{equation}\label{locus_sing}
  {b}^{2} \left( -{b}^{2}+3\,{a}^{2} \right) ^{2} \left( {a}^{2}+{b}^{2}  \right) ^{6}=0
\end{equation}
  which leads us to the following analysis. It is sufficient to study solutions of (\ref{locus_sing}) on the unit circle $a^2+b^2=1$, so we must have
   $$( a-1) ( a+1) ( 2\,a-1) ^{2}( 2\,a+1) ^{2}=0.$$
Studying the various possibilities we find the following.
 \begin{enumerate}\itemindent=4pt\itemsep=4pt
   \item[(a)] There is a $1$-cell of triple self-intersection points of $Z(\Delta_3)$ corresponding to  polynomials of form $z^3+2\,d\,i$,
       with $d$ varying in $\R$ (see Figure~\ref{fig_forest_trois}-(a)).  One calculates that the matching characterizing this cell is 
     $\{\bleu{\{1,7\},\{3,9\},\{5,11\}}\}$.
   \item[(b)] Moreover there are three $2$-cells  decomposing the rest of the locus of self-intersections of the surface. These correspond respectively to  polynomials of the form
\begin{eqnarray*}
     z^3+3\,r\,z-2\,d\,i, &&  z^3-3\,re^{i\,\pi/3}\,z+2\,d\,i,\quad \mathrm{and}\\
   && z^3-3\,re^{-i\,\pi/3}\,z+2\,d\,i,
\end{eqnarray*}
   where we assume that $r>0$, and $d\in \R$. See Figure~\ref{fig_forest_trois}-(b) for a polynomial that is of the first form. The corresponding matchings are 
$$\{\bleu{\{1,5\},\{3,9\},\{7,11\}}\},\quad
 \{\bleu{\{1,7\},\{3,11\},\{5,9\}}\},\quad \mathrm{and}\quad 
  \{\bleu{\{1,9\},\{3,7\},\{5,11\}}\}.$$
      \item[(c)] Outside of this locus of self-intersection, there are six  $3$-cells consisting of polynomials of the form
         $$z^3-3\,r\,e^{3\,i\,\theta/2}\,z+2 (r^{3/2}\cos(3\,\theta/2) +d\,i),$$
again with $r>0$, $d\in\R$. Here we assume that $\theta$ satisfies lies in one of the open intervals
     $$ \begin{array}{ll}
          \mathrm{(1)}\    -\pi/3<\theta<\pi/3 ,\qquad & \mathrm{(2)}\  \pi/3<\theta<\pi ,\\[4pt]
          \mathrm{(3)}\    \pi<\theta<5\,\pi/3 ,\qquad & \mathrm{(4)}\   5\,\pi/3<\theta<7\,\pi/3 ,\\[4pt]
          \mathrm{(5)}\    7\,\pi/3<\theta<3\,\pi ,\qquad & \mathrm{(6)}\  3\,\pi/3<\theta<11\,\pi/3, 
       \end{array}$$
 which correspond respectively to one of the $3$-cells. The first interval discribes a $3$-cell in which lies the polynomial $z^3-3\,z+2$, obtained by choosing $\theta=0$, $r=1$, and $d=0$. The corresponding forest appears in Figure~\ref{fig_forest_trois}-(c).
 The six matchings that correspond to these cells are 
$$\begin{array}{rrrrr}
\{\bleu{\{1,9\},\{3,11\},\{5,7\}}\}, & \{\bleu{\{1,5\},\{3,11\},\{7,9\}}\},\\[4pt]
\{\bleu{\{1,5\},\{3,7\},\{9,11\}}\}, & \{\bleu{\{1,11\},\{3,7\},\{5,9\}}\},\\[4pt]
\{\bleu{\{1,3\},\{5,9\},\{7,11\}}\}, & \{\bleu{\{1,9\},\{7,11\},\{3,5\}}\}.
\end{array}$$
     \end{enumerate} 
We conclude that the resulting cell enumerator is
\begin{equation}
   \mathcal{F}_3= x+3\,x^2+6\,x^3+5\,x^4.
\end{equation}
We also deduce from our discussion that there are $5$ cell orbits, and that their degree weighted enumeration gives
  $$x+x^2+x^3+ 2\,x^4,$$

\section{Degree 3 polynomials: basketballs}\label{cas_bask_3}
\begin{figure}[ht]
 \begin{center}
\includegraphics[width=40mm]{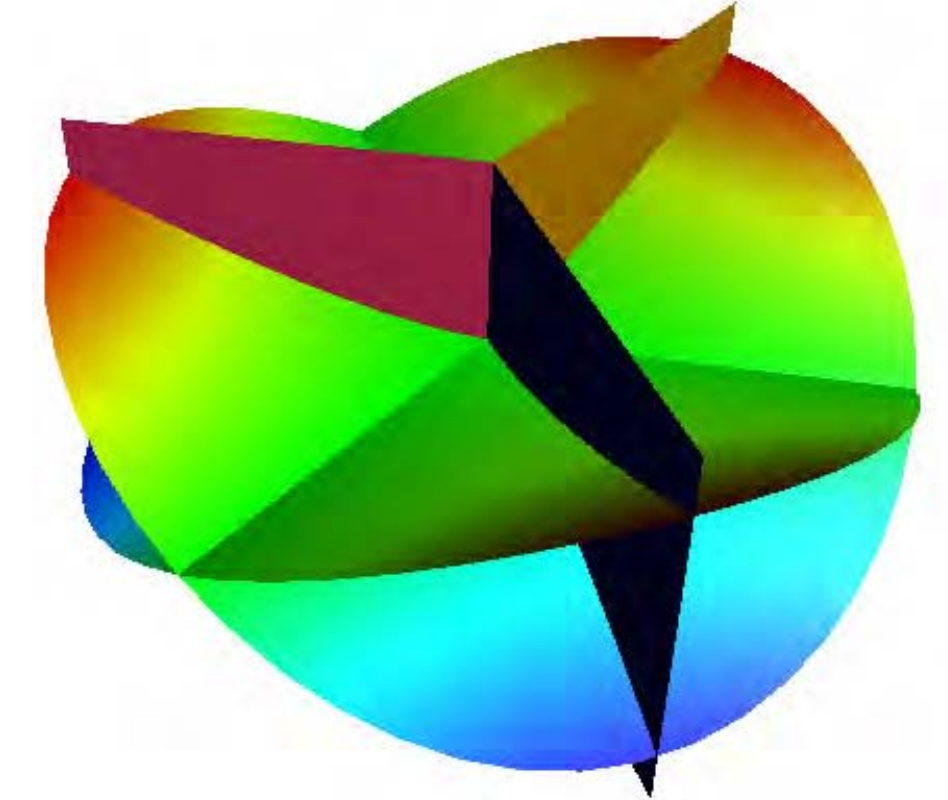}
\end{center}\vskip-10pt
\caption{Basketball singularity locus at $c=0$.}\label{fig_d=0}
 \end{figure}

The next order of buisness is to further decompose the space $\mathcal{P}_3$ taking into account the pairs $(\Rf(f),\imf(f))$. Thus we break up the cells of the previous decomposition by considering the  imaginary component. We are going to organize our presentation around this further decomposition.
To start this process, one calculates that in cylindrical coordinates
 \begin{equation}\label{delta_prime}
      \Delta_3'=4\,(d^2 +r^3\cos^2( 3\,\theta/2) ) (d^2-r^3\sin^2( 3\,\theta/2)), 
  \end{equation} 
 which expression is now independent of $c$. Again, since $r>0$ and $d$ is real, so that we have $\Delta'\geq 0$ if and only if
 \begin{equation}\label{surface_prime}
     d^2\geq r^3\sin^2( 3\,\theta/2).
  \end{equation}
To summarize, the singularity locus $Z(\Delta_3\Delta_3')$ for basketballs of 
   $$f(z)=z^3-3\,(a+b\,i)z+2\,(c+d\,i)$$
is the set of polynomials either of the form
         $$z^3-3\,r\,e^{i\,\theta}\,z+2 (r^{3/2}\cos(3\,\theta/2) +d\,i),$$
or
         $$z^3-3\,r\,e^{i\,\theta}\,z+2 (c+r^{3/2}\sin(3\,\theta/2)\,i),$$
 with $r>0$, and $0\leq \theta \leq 4\,\pi$. Restricting the singularity locus to the hyperplane $c=0$, we get the surface illustrated in Figure~\ref{fig_d=0}.
This corresponds to restricting the real component of the basketball to be either of type (b), (c) or (d) in  our classification of the real component of section~\ref{sec_real}.  
 
We already know that there are $22$ cells corresponding to non-singular basketballs, and these are the only $4$-cells. 
Recall that these cells come in a group of three orbits, that are respectively generated by families of polynomials of the form 
$$z^3-3\,r\,e^{i\,\theta}\,z+2 r^{3/2}(c +d\,i),$$
for $r>0$. We assume that $\theta$ is chosen so that both  $\cos(3\,\theta/2)>0$ and $\sin(3\,\theta/2)>0$. Then the generating cells are
\begin{enumerate}\itemindent=4pt\itemsep=4pt
  \item[(A.1)] A cell for which $c> \cos(3\,\theta/2)$,  and $d>\sin(3\,\theta/2)$ (see Figure~\ref{vingtdeux}-(A.1)). This $4$-cell generates an orbit of length $4$. The associate matching is 
       $$\{\bleu{\{ 1,3\} ,\{ 5,7\} ,\{ 9,11\} },\rouge{\{ 0,2\} ,\{ 4,6\} ,\{ 8,10\}}\}.$$
  \item[(A.2)] A cell for which $c> \cos(3\,\theta/2)$,  and $-\sin(3\,\theta/2)<d<\sin(3\,\theta/2)$ (see Figure~\ref{vingtdeux}-(A.2)). The generated orbit has length $12$, and the matching that corresponds to the generating cell is
        $$\{\bleu{\{ 1,3\} ,\{ 5,7\} ,\{ 9,11\} },\rouge{\{ 0,2\} ,\{ 4,10\} ,\{ 6,8\}}\}.$$
    \item[(B)]  A cell for which $- \cos(3\,\theta/2)<c<\cos(3\,\theta/2)$, and $-\sin(3\,\theta/2)<d<\sin(3\,\theta/2)$ (see Figure~\ref{vingtdeux}-(A.2)). The corresponding orbit has length $6$. The matching characterizing the generating cell is now
        $$\{\bleu{\{ 1,11\} ,\{ 3,9\} ,\{ 5,7\} },\rouge{\{ 0,2\} ,\{ 4,10\} ,\{ 6,8\}}\}.$$    
 \end{enumerate}
 The remaining  cells correspond to singular basketballs, and the corresponding decomposition naturally breaks up into orbits under the action of the relevant dihedral group.  Recall that this dihedral group is generated by the rotation $f(z)\mapsto i\,f(e^{-\pi i/2n}\,z)$ that exchanges the role of $\Rf$ and $\imf$, and the reflection $\overline{ f(\overline{z})}$. We classify these orbits using representatives of the forest-cell-decomposition, which are further decomposed to get representative of each of the possible orbits of basketballs.
The orbits are as follows:
\begin{enumerate}\itemindent=4pt\itemsep=4pt
  \item[(a.1)] Inside the $1$-cell corresponding to case (a) (of the real component forest-cell-decomposition) sits the $0$-cell consisting of the unique polynomial $z^3$. This is the only $0$-cell of the whole basketball-cell-decomposition, and its orbit is trivial. 
    \item[(a.2)] Together with this $0$-cell, two more $1$-cells (corresponding to cases for which $a=b=0$ and $c\,d=0$) funrnish the entire basketball-cell-decomposition of the forest-cell (a). There is just one length $4$ orbit generated by any of these two $1$-cells. These cells are respectively made up of  the polynomials 
       $$ \begin{array}{llll}
              z^3-2\,d\,i,&  z^3+2\,d\,i, &   z^3-2\,c,\ \mathrm{and}  & z^3+2\,c,
           \end{array}$$
for $c>0$ and $d>0$.   To summarize, we have three cells $z^3+2\,d\,i$, $z^3$, and $z^3-2\,d\,i$ (with $d>0$)
decomposing (a), and they belong to one of the two orbits generated respectively by
 \begin{center}
\includegraphics[width=50mm]{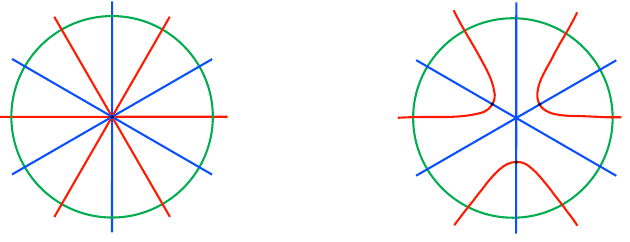}
\begin{picture}(0,0)(150,10)\small
\put(0,0){(a.1)\  $z^3$} 
\put(90,0){(a.2)\   $z^3-2\,d\,i$}  
\end{picture}\vskip15pt
\end{center}
These give the respective matchings
\begin{eqnarray*}
&& \{\bleu{  \{ 1,7 \} , \{ 3,9 \} , \{ 5,11 \} } ,\rouge{ \{ 0,2 \} , \{ 4,6 \} , \{ 8,10 \}}\},\quad \mathrm{and}\\
&&\{ \bleu{ \{ 1,5 \} , \{ 3,9 \} , \{ 7,11 \} } ,\rouge{ \{ 0,4 \} , \{ 2,6 \} , \{ 8,10 \}}\}.
\end{eqnarray*}
 
   \item[(b.1)] Recall that all of the two forest-cells in (b) can be obtained from one another by the action of the group. Thus, as our next step, we need only study how to break up the first one of these. This is the one corresponding to  polynomials of the form
           \begin{equation}\label{cas_b}
                  z^3+3\,r^2\,z-2\,d\,i.
          \end{equation}
To this aim, we first consider the orbit generated by the $1$-cell, sitting inside (b), consisting of polynomials  that have a double root. These are polynomials in the orbit of
   $$z^3+3\,r^2\,z-2\,i\,r^3,$$
 for $r>0$.   This cell generates the following orbit of $12$ cells:
   $$\begin{array}{lll} 
  {z}^{3}+3\,r^2 z-2\,i\,r^3,\qquad&
  {z}^{3}+3\,q\,r^2 z+2\,r^3, \qquad&
  {z}^{3}-3\,\overline{q}\,r^2 z+2\,i\,r^3,\\[4pt]
  {z}^{3}+3\,r^2 z+2\,i\,r^3,&
  {z}^{3}+3\,q\,r^2 z-2\,r^3,&
  {z}^{3}-3\,\overline{q}\,r^2 z-2\,i\,r^3,\\[4pt]
  {z}^{3}-3\,r^2z+2\,r^3,  &
  {z}^{3}-3\,q\,r^2\, z+2\,i\,r^3, &
  {z}^{3}+ 3\,\overline{q}\,r^2\, z-2\,r^3,\\ [4pt]
  {z}^{3}-3\,r^2 z-2\,r^3,&
  {z}^{3}-3\,q\,r^2\, z-2\,i\,r^3,&
  {z}^{3}+ 3\,\overline{q}\,r^2\, z+2\,r^3,
  \end {array}$$
      with $q$ denoting the third root of unity $\exp(i\,\pi/3)$. Only the first two sets of the first column actually contribute to the decomposition of the specific forest-cell under consideration. 
   \item[(b.2)] To complete our decomposition of the cell containing the polynomials in  (\ref{cas_b}),   we next consider those for which
 with $d>r^3$, so that we now have a $2$-cell. Again, we get an orbit of length $12$. Once more, inside this orbit, only two cells do contribute to the decomposition of the cell. Indeed, beside those already mentioned we also need to consider polynomials of form  (\ref{cas_b}) such that
 $d<-r^3$.
  
    \item[(b.3)] To finish our decomposition of the cell corresponding to polynomials of form (\ref{cas_b}), it only remains to consider the orbit of polynomials for which $-r^3<d<r^3$. The corresponding cell is also a $2$-cell, and it now generates an orbit of length $6$. This is made apparent by considering the corresponding basketball below. Altogether, we have considered the orbits of the three  cells.
    \begin{center}
\includegraphics[width=100mm]{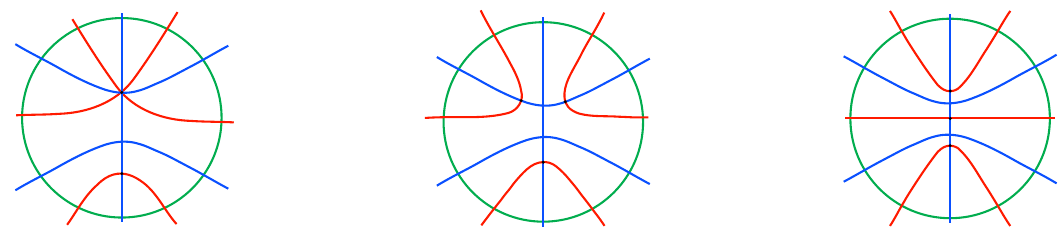}
\begin{picture}(0,0)(320,10)\small
\put(0,0){(b.1)\  $z^3+3\,r^2\,z-2\,i\,r^3$} 
\put(120,0){(b.2)\   $z^3+3\,r^2\,z-2\,i\,d$}  
\put(150,-15){($d>r^3$)}  
\put(240,0){(b.3)\   $z^3+3\,r^2\,z-2\,i\,d$} 
\put(270,-15){($-r^3<d<r^3$)}  
\end{picture}\vskip25pt
\end{center}
The matchings that correspond to these are now
\begin{eqnarray*}
&&\{ \bleu{ \{ 1,5 \} , \{ 3,9 \} , \{ 7,11 \} } ,\rouge{ \{ 0,4 \} , \{ 2,6 \} , \{ 8,10 \}}\},\\
&& \{ \bleu{ \{ 1,5 \} , \{ 3,9 \} , \{ 7,11 \} } ,\rouge{ \{ 0,2 \} , \{ 4,6 \} , \{ 8,10 \}}\},\\
&&\{ \bleu{ \{ 1,5 \} , \{ 3,9 \} , \{ 7,11 \} } ,\rouge{ \{ 0,6 \} , \{ 2,4 \} , \{ 8,10 \}}\}.
\end{eqnarray*}
   \item[(c.1)] We are only left with the analysis of how to decompose the forest cells classified in case (c), and we do this for  
    the forest-cell containing the polynomials
\begin{equation}\label{cas_c}
   z^3-3\,r\,e^{i\,\theta}\,z+2\,r^{3/2}(e^{3\,i\,\theta/2} +  d\,i),
 \end{equation}
with $r>0$ and $0<\theta<\pi/3$, $d\in\R$. As before, we decompose this cell by studying the different basketballs that can occur as $d$ varies from $-\infty$ to $+\infty$. The first sub-cell encountered generates a length $24$ orbit of $3$-cells, and it contains the polynomials of form (\ref{cas_c}) for which $d<-2\sin(3\,\theta/2)$.
   \item[(c.2)]  Then one reaches one orbit  of $2$-cells, of length $12$, which contains polynomials for which both the real component and the
   imaginary component exhibit a singularity. This occurs exactly when $d=-2\sin(3\,\theta/2)$. 
   \item[(c.3)]  After that, for $-2\sin(3\,\theta/2) <d<0$, we get a $3$-cell whose orbit is of length $24$,

   \item[(c.4)]  and then we reach $d=0$, at which point the polynomials under study have a double root. The orbit of the corresponding $2$-cell is of length $12$. The last cell in the decomposition of case (c) corresponds to  $d>0$,  but this cell has already been taken care of since it appears in  the orbit considered in (c.1). Thus the whole situation corresponds to cells associated to an orbit of one of the following cells:
    \begin{center}
\includegraphics[width=70mm]{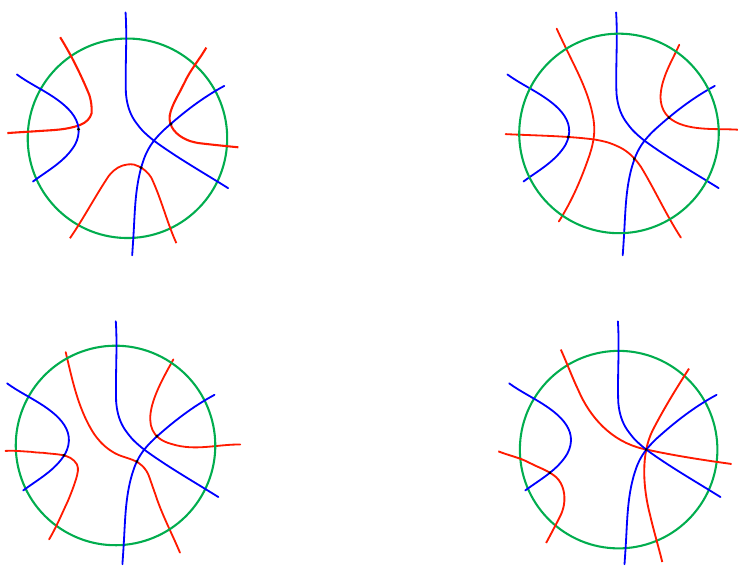}
\begin{picture}(0,0)(350,25)\small
\put(135,100){(c.1)\ $d<-2\sin(3\,\theta/2)$}   
\put(260,100){(c.2)\  $d=-2\sin(3\,\theta/2)$}  
\put(130,15){(c.3)\  $-2\sin(3\,\theta/2) <d<0.$} 
\put(285,15){(c.4)\  $d=0.$} 
\end{picture}\vskip15pt
\end{center}
and the corresponding matchings are
\begin{eqnarray*}
&&\{\bleu{  \{ 1,9 \} , \{ 3,11 \} , \{ 5,7 \} } ,\rouge{ \{ 0,2 \} , \{ 4,6 \} , \{ 8,10 \}}\},\\
&&\{\bleu{  \{ 1,9 \} , \{ 3,11 \} , \{ 5,7 \} } ,\rouge{ \{ 0,2 \} , \{ 4,8 \} , \{ 6,10 \}}\},\\
&&\{\bleu{  \{ 1,9\} , \{ 3,11 \} , \{ 5,7 \} } ,\rouge{ \{ 0,2 \} , \{ 4,10 \} , \{ 6,8 \}}\},\\
&&\{\bleu{ \{ 1,9 \} , \{ 3,11 \} , \{ 5,7 \} } ,\rouge{ \{ 0,4 \} , \{ 2,10 \} , \{ 6,8 \}} \}.
\end{eqnarray*}
    \end{enumerate}
Adding all this up, including the $22$ non-singular basketballs cells of degree $4$, we get the cell enumerator:
  $$\mathcal{B}_3=1+16\,x+42\,{x}^{2}+48\,{x}^{3}+22\,{x}^{4}.$$
Under the action of the subgroup that respects the basketball structure (without exchanging real and imaginary component), the degree enumeration of orbits results in the polynomial
 $$1+4\,x+8\,{x}^{2}+8\,{x}^{3}+6\,{x}^{4}.$$

\section{Polynomials of the form $z^n-n\,(a+b\,i)z+(n-1)\,(c+d\,i)$}  \label{borne_un}
More generally, restricting our discussion to the subspace $\mathcal{Q}_n$ of polynomials   of the form 
      $$f(z)=z^n-n\,(a+b\,i)z+(n-1)\,(c+d\,i)$$ 
 we find that the singularity locus for singular $\Rf(f)$ is the set
$$Z(\Delta_n)=\left\{( r\cos( \theta) ,r\sin( \theta) , r ^{{n}/{(n-1)}}\cos
( n\,\theta/(n-1)),d)\ |\ r\geq 0,\ 0\leq \theta< 2\,(n-1)\,\pi\right\},$$
in cylindrical coordinates.
Once again we observe that there is no condition on the $d$-coordinate, and thus the polynomials that exhibit singularity in the real component of their basketballs are  those that can be written as
\begin{equation}\label{singular_general}
   f(z)=z^n-n\,r\,e^{i\,\theta} z+(n-1)\,r^{n/(n-1)}\cos(n\,\theta/(n-1))+d\,i.
\end{equation}
This fact can made apparent by the verification that the critical values of $z^n-n\,r\,e^{i\,\theta} z$ are all of the form $(n-1)\,r^{n/(n-1)}e^{\,n\,i\, \theta/(n-1)}$,
with $\theta$ specified modulo $2\,\pi$. Since $d$ may take arbitrary values without changing the real component of the basketball, we may suppose that it ``absorbs''  the imaginary part of the critical value. The real parts of the critical values 
   $$c_k:=(n-1)\,r^{n/(n-1)}\cos\left(\frac{n}{n-1}\,\theta+\frac{2\,k}{n-1}\,\pi\right),\qquad k=0,\ldots,n-2,$$
are all distinct for generic values of $\theta$. However, when $\theta=j\,\pi/n$, some of the $c_k$ will clearly coincide, and  the surface $Z(\Delta_n)$ will self-intersect.

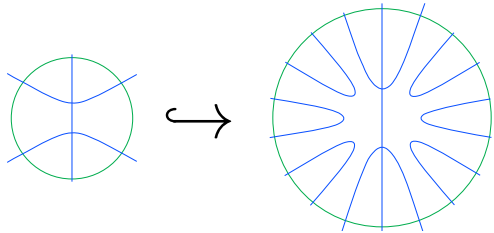
\begin{figure}[ht]
 \begin{center}
\includegraphics[width=70mm]{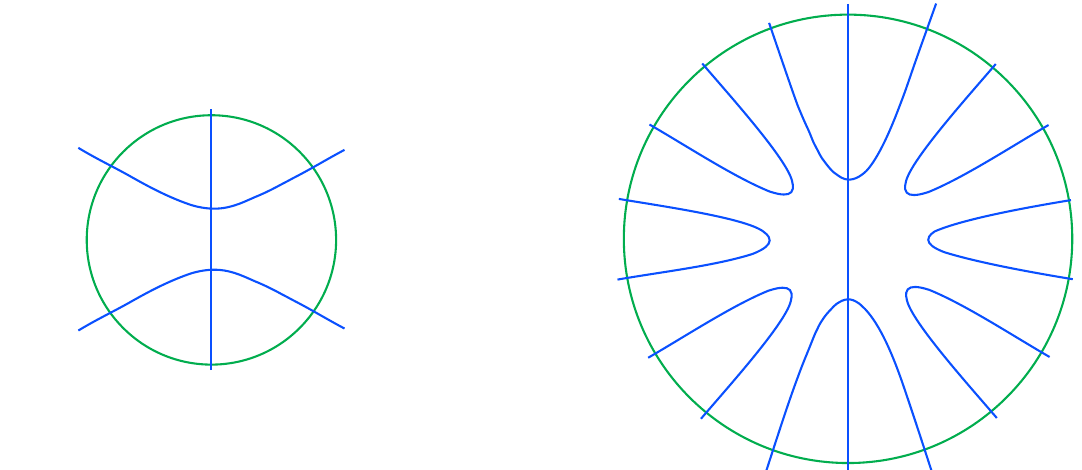}
\begin{picture}(0,0)(0,0)
\put(-130,35){\Huge $\hookrightarrow$}
\end{picture}
\end{center}\vskip-10pt
\caption{A $9$-embedding of a $5$-matching.}\label{fig_ n_general}
 \end{figure}

It is easy to verify that the restriction to $\R^3$ of $Z(\Delta_n)$ gives a surface that decomposes $\R^3$ in $n\,(n-2)+2$ open cells, $n\,(n+1)$ cells of degree $2$, $n\,(n-2)$ cells of degree $1$, and one $0$-cell, for a total of $3\,n^2-5\,n+3$ cells. We simply add $1$ to the degree of these cells  to get the forest cell enumerator of $\mathcal{Q}_n$ to get
\begin{equation}\label{enumerator_n_2}
   x+(n^2-2\,n)\,x^2+(2\,n^2-4\,n)\,x^3+(n^2-2\,n +2)\,x^4.
 \end{equation}
 To describe the possible circular forests that correspond to polynomials in the space $\mathcal{Q}_n$, we consider the notion of $n$-embedding of circular $k$-trees (with $2\,k$ leaves), for $k=1,2,3$. For $n>k$, a $n$-{\em embedding} of a $k$-tree $T$  is an increasing function 
    $$\varphi:\{1,3,\ldots,4\,k-1\} \rightarrow \{1,3,\ldots,4\,n-1\},$$
such that there is a unique circular forest $F$ on the $2\,n$ points of the bigger set with the following conditions holding
\begin{enumerate}\itemindent=4pt\itemsep=4pt
\item[(1)] the branch $\{\varphi(a),\varphi(b)\}$ is in $F$ for all branch $\{a,b\}$ is in $T$, and 
\item[(2)] all other branches are ``trivial''. This is to say that for any two consecutive\footnote{There are no other points in the image of $\varphi$ that lie between $c$ and $e$ going counterclockwise.} values $c$ and $e$ in the image of $\varphi$, if $\{d_k\}$ is in the set of
 points that lie in between $c$ and $e$, ordered counterclockwise, then $d_k$ is matched to $d_{k+2}$.
 \item[(3)] If $k=2$  at most one leaf in the image is not consecutive to the others. If $k=3$, all leaves attached to the same internal node are sent to consecutive leaves by $\varphi$.
\end{enumerate}
An example of a $9$-embedding of a $3$-tree is given in Figure~\ref{fig_ n_general}.
It follows from our discussion that
  \begin{prop}
Except for the real component of basketballs corresponding to polynomials of form   $z^n+ d\,i$,
   the  real component of basketballs of polynomials in $\mathcal{Q}_n$ are all obtained as $n$-embeddings of circular $k$-trees {\rm (}with $k=1,2,3${\rm )}  whose internal vertices have at most $4$ neighbours. Calculating the forest cell enumerator of the space $\mathcal{Q}_n$ using this characterization gives back formula {\rm (}\ref{enumerator_n_2}{\rm )}.
\end{prop}

 A similar study can be done for basketballs, and all the connected basketballs appearing in a basketball describing the relevant cells already appear for $n\leq 5$.

\section{Circular forest enumeration}
\begin{figure}[ht]
 \begin{center}
\includegraphics[width=40mm]{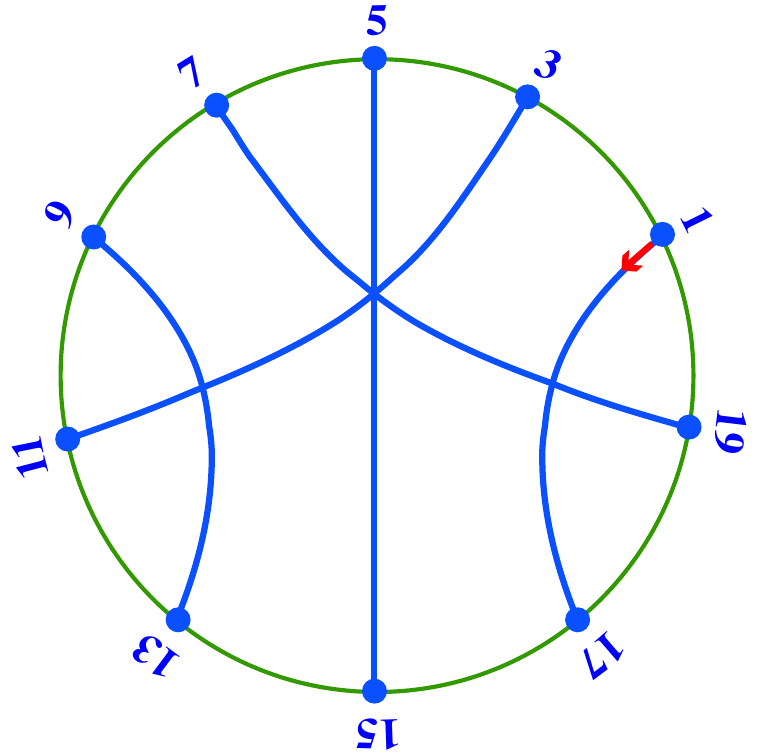}
\end{center}\vskip-10pt
\caption{Rooted circular tree.}\label{fig_arbre}
 \end{figure}

We are now going to enumerate general circular forests. As a first step, we consider the enumeration of circular trees with even degrees for internal nodes. Observe that such  trees can be canonically ``rooted'' at $1$ (see Figure~\ref{fig_arbre}).
Orienting edges (to get arcs), going up from the root toward the leaves, each internal vertex has one ingoing arc and an odd number ($\geq 3$) of outgoing arcs.  Now, consider the generating function of the series
     $$\mathcal{A}(x,t):=\sum_{n} a_{n}(t)\, x^{2\,n-1},\qquad \mathrm{with}\qquad a_n(t)=\sum_{k} a_{n,k}\,t^k$$
 and $a_{n,k}$ equal to the number of rooted circular trees that have $k$ internal nodes and $2\,n-1$ leaves (not counting the root). It follows from the above observation that $\mathcal{A}=\mathcal{A}(x,t)$ satisfies the functional equation (see \cite{BLL} for further details on how to establish this)
\begin{equation}\label{funct_eq}
    \mathcal{A}=x + \frac{t\,\mathcal{A}^3}{1-\mathcal{A}^2}.
\end{equation}
From this, one calculates the series expansion
\begin{eqnarray*}
\mathcal{A}(x,t)&=&x+t\,{x}^{3}+ ( t+3\,{t}^{2})\, {x}^{5}+
 ( t+8\,{t}^{2}+ 12\,{t}^{3}) \,{x}^{7}+
  ( t+15\,{t}^{2}+55\,{t}^{3}+55\,{t}^{4} )\, {x}^{9}\\
&&\qquad + (t+24\,{t}^{2}+ 156\,{t}^{3}+364\,{t}^{4}+273\,{t}^{5} )\, {x}^{11}\\
&&\qquad + (t+35\,{t}^{2}+350\,{t}^{3}+ 1400\,{t}^{4}+2380\,{t}^{5}+1428\,{t}^{6} )\, {x}^{13}+\ldots
\end{eqnarray*}
It is clear that a circular forests consists in a choice of an even parts non-crossing partition of the $2\,n$ points on the circle, and a choice of a circular tree
structure on each of the parts of this partition.  Thus if we count non-crossing partitions $\lambda$ of $2\,n$ points with a  weight equal to 
     $$\omega(\lambda)=x^{2n}\,\prod_{B\in\lambda} p_{\#{B}} ,$$
then we can enumerate circular forest  by the simple process of replacing each $b_{2k}$ by $a_k=a_k(t)$.
Now, the weighted enumeration of even parts non-crossing partitions results in the series
\begin{eqnarray*}
\Theta(x)&=&  {p_2}\,x^2+( 2\,p_2^{2}+{p_4})\, {x}^{4}+( 5
\,p_2^{3}+6\,{p_2}\,{p_4}+{p_6})\, {x}^{6}\\
&&+
( 14\,p_2^{4}+28\,p_2^{2}\,{p_4}+8\,{p_2}\,{p_6}+4\,p_4^{2}+{p_8})\, {x}^{8}\\
&& +( 42\,p_2^{5}+120\,p_2^{3}\,{p_4}+45\,p_2^{2}\,{p_6}+45\,
 {p_2}\,p_4^{2}+10\,{p_2}\,{p_8}+10\,{p_4}\,{p_6}+
{p_{10}})\, {x}^{10}+\ldots
\end{eqnarray*}
It is easy to recursively generate representatives of non-crossing partitions with respect to the action of the group of rotations. Indeed there is always one part of the partition consisting of consecutive points on the circle. Up to a rotation, we can assume that this part consists in the $k$ last points in our labelling. Removing this part, we are thus left with a smaller non-crossing partition.
Using the above series expansion, we find that the generating function for circular forest is
\begin{eqnarray*}
\mathcal{F}(x,t)&=&x^2+ \left( 2+t \right) {x}^{4}+ \left( 5+7\,t+3\,{t}^{2} \right) {x}^{6
}+ \left( 14+37\,t+36\,{t}^{2}+12\,{t}^{3} \right) {x}^{8}\\
&&\qquad + \left( 42+
176\,t+285\,{t}^{2}+205\,{t}^{3}+55\,{t}^{4} \right) {x}^{10}\\
&&\qquad + \left( 
132+794\,t+1872\,{t}^{2}+2158\,{t}^{3}+1222\,{t}^{4}+273\,{t}^{5}
 \right) {x}^{12}+\ldots
  \end{eqnarray*}
where the power of $t$ marks the number of internal nodes, and that of $x$ marks the number of leaves.  

\section{Further considerations}
For all fixed $k$ and $n>K$, one can extend to polynomials of the form 
   $$z^n+\alpha_k\, z^k+\ldots + \alpha_1\,z+\alpha_0,$$
 many of the considerations of section~\ref{borne_un}. In particular, the non-crossing matchings that are associated to these polynomials are such that
 at most $k$ pairs consist in two points that are not immediate neighbours. It follows that all possible crossing matchings that arise for $n$ can be obtained as embeddings of  crossing matchings belonging to  a fixed finite list depending only on $k$. Many other special families of polynomials (such as families of orthogonal polynomials, or conservative polynomials (see \cite{smale}))  have basketballs  exhibiting nice combinatorial properties. To illustrate what we mean here, consider the family of real-coefficient polynomials $f(z)$ having distinct real roots, so that the roots of $f'(z)$ are interlaced with the roots of $f(z)$. This family corresponds exactly to one of the cells of the basketball  cellular decomposition.  This cell is characterized by the matching for which $0$ is paired with $2\,n$, and any other $k$ is paired with $-k$ (modulo $4\,n$). From what we have already said, it is clear that the intersection of the branches $\{k,-k\}$ (modulo $4\,n$) with the branch $\{0,2\,n\}$ is a zero of $f(z)$ if $k$ is odd, and a zero of $f'(z)$ if $k$ is even. Hence the real component is noncrossing, and the imaginary component is a tree. See \cite{lasserre} on how to get a semi-algebraic characterization of this cell in terms of principal minors of Hankel matrices with entries equal to sums of power of roots.
  
One of the remaining problem is that of effectively constructing polynomials realizing given combinatorial basketballs.  Even better, one would like to have an explicit semi-algebraic description of cells associated to their combinatorial basketball indexing. 
Along these lines, a recent theorem of \cite[see Thm. 1.3]{jakobson} essentially  states that it is always possible to find a polynomial realizing a given circular forest (as the real component of its basketball). Hence this can be considered as a partial generalization of the ``Inverse Basketball Theorem'' of \cite{martin} with regards to real components. We expect that one should be able to derive that all basketballs are indeed realizable  from the techniques used in \cite{jakobson} and \cite{martin}.  On the other hand, the question of effectively giving semi-algebraic descriptions of cells is not being addressed by these results.


\end{document}